\documentclass[a4paper]{amsart}

\DeclareMathOperator{\Aut}{Aut} \DeclareMathOperator{\spn}{span}

\newcommand{\Z}{\mathbb{Z}}
\newcommand{\N}{\mathbb{N}}
\newcommand{\T}{\mathbb{T}}
\newcommand{\mc}[1]{\mathcal{#1}}

\newcommand{\map}[3]{#1:#2\rightarrow#3}
\newcommand{\set}[1]{\left\{#1\right\}}
\newcommand{\OR}{\mbox{ or }}
\newcommand{\cspan}{\overline{\mathrm{span}}} 
\newcommand{\ol}{\overline}

\newcommand{\CSalgebra}[1][\ ]{$C^*$-algebra#1}
\newcommand{\CSalgebras}[1][\ ]{\CSalgebra[s]#1}
\newcommand{\CS}[1]{C^*\left(#1\right)}

\newcommand{\CKLXfam}[1][\ ]{Cuntz-Krieger $(\Lambda,X)$-family#1}
\newcommand{\CKLXfams}[1][\ ]{Cuntz-Krieger $(\Lambda,X)$-families#1}

\newcommand{\KGalg}[1][\Lambda]{\CS{#1}}
\newcommand{\KGalgC}[2][\Lambda]{\CS{#1,#2}}
\newcommand{\KGalgCLX}{\KGalgC{X}}
\newcommand{\KGalgCLL}{\KGalgC{\Lambda^0}}

\newcommand{\XsubLam}[1][0]{X\subseteq\Lambda^{#1}}

\newcommand{\Corner}[2]{P_{#2}\CS{#1}P_{#2}}
\newcommand{\CorLX}{\Corner{\Lambda}{X}}

\newcommand{\Lambdamin}[3][\Lambda]{#1^{\min}(#2,#3)}

\newcommand{\Tl}[2][T]{#1_{#2}}
\newcommand{\Tlam}[1][T]{\Tl[#1]{\lambda}}
\newcommand{\Tab}[3][T]{\Tl[#1]{#2,#3}}
\newcommand{\Talbe}[1][T]{\Tab[#1]{\alpha}{\beta}}

\newcommand{\sat}[1]{\Sigma\left(#1\right)}

\newcommand{\SPg}[3][c]{#2\times_{#1}#3}
\newcommand{\SPgGL}{\SPg{G}{\Lambda}}
\newcommand{\SPgZL}{\SPg[d]{\Z^k}{\Lambda}}
\newcommand{\KGalgCGLV}{\KGalgC[\SPgGL]{V}}
\newcommand{\KGalgCZLV}{\KGalgC[\SPgZL]{V}}

\newcommand{\St}[1]{\mc{S}\left(#1\right)}
\newcommand{\StL}{\St{\Lambda}}
\newcommand{\ess}[1]{E\left(#1\right)}
\newcommand{\eLambda}{\ess{\Lambda}}

\theoremstyle{plain}
\newtheorem{thm}{Theorem}[section]
\newtheorem{prop}[thm]{Proposition}
\newtheorem{lem}[thm]{Lemma}
\newtheorem{cor}[thm]{Corollary}

\theoremstyle{definition}
\newtheorem{defn}[thm]{Definition}

\theoremstyle{remark}
\newtheorem{rem}[thm]{Remark}
\newtheorem{remarks}[thm]{Remarks}
\newtheorem{eg}[thm]{Example}
\newenvironment{rems}{\begin{remarks}\begin{enumerate}}{\end{enumerate}\end{remarks}}

\newcommand{\textref}[3][\ref]{#2\ #1{#3}}
\newcommand{\defref}[1]{\textref{Definition}{#1}}
\newcommand{\thmref}[1]{\textref{Theorem}{#1}}
\newcommand{\propref}[1]{\textref{Proposition}{#1}}
\newcommand{\lemref}[1]{\textref{Lemma}{#1}}
\newcommand{\corref}[1]{\textref{Corollary}{#1}}
\newcommand{\egref}[2][\ref]{\textref[#1]{Example}{#2}}

\usepackage{amssymb}

\begin{document}

\title{Gauge Invariant Uniqueness Theorem for Corners of $k$-graphs}
\author{Stephen Allen}
\date{\today}
\email{stephen.allen@newcastle.edu.au}
\address{School Mathematical and Physical Sciences, University of Newcastle, Callaghan, NSW 2308, AUSTRALIA}
\subjclass[2000]{46L05}

\begin{abstract}
For a finitely aligned $k$-graph $\Lambda$ with $X$ a set of
vertices in $\Lambda$ we define a universal \CSalgebra called
$\KGalgCLX$ generated by partial isometries. We show that
$\KGalgCLX$ is isomorphic to the corner $\CorLX$, where $P_X$ is
the sum of vertex projections in $X$. We then prove a version of
the Gauge Invariant Uniqueness Theorem for $\KGalgCLX$, and then
use the theorem to prove various results involving fullness,
simplicity and Morita equivalence as well as new results involving
symbolic dynamics.
\end{abstract}

\maketitle

\section{Introduction}
Much study has been done lately in regards to higher rank graphs
(also known as $k$-graphs) and their associated graph algebras
since their first appearance in \cite{KP2}. As $k$-graphs are a
higher dimensional generalisation of directed graphs (which can be
assumed to be 1-dimensional), it is important to be able to adapt
the known results for directed graphs to the field of $k$-graphs.
So far this has been done with a reasonable amount of success (for
example see \cite{APS} \cite{RSY1}, \cite{KP3}, \cite{Si2} to name
a few) however the complex nature of of $k$-graphs often makes the
proofs of these adapted results much more complicated than the
previous ones.

Corners of graph algebras naturally arise in many places when
studying graph algebras (see \cite{DT2}, \cite{KPRR}, \cite{Ty},
\cite{Tom1} for example) and have shown to be a necessary tool in
the understanding of arbitrary graph algebras. In particular,
there is an important link between graph algebras and symbolic
dynamics (see \cite{Ba}, \cite{BP} and \cite{DS}) since directed
graphs represent subshifts of finite type (see \cite{LM}).
Transferring results from symbolic dynamics to graph algebras
frequently involves using corners.

It is the goal of this paper to provide tools for dealing with
corners of $k$-graph algebras generated by vertex projections. As
such we describe a universal \CSalgebra generated by partial
isometries which is isomorphic to a corner of a graph algebra. We
then show that this algebra has a version of the Gauge Invariant
Uniqueness Theorem (\thmref{thm:GUIT}) which tells us when
mappings that respect the gauge action are injective. We then show
the facility our definition provides by proving various
applications. As such we gain some new results as well as some
generalisations of existing results for directed graphs. In
particular we obtain conditions for checking Morita equivalence of
graph algebras using corners and also realise the AF core of a
$k$-graph as a corner.

We begin in section 2 with the preliminaries involved in finitely
aligned $k$-graphs and their associated graph algebra since for
the most part of this paper we restrict ourselves to this class of
$k$-graph. In section 3, given a set $X$ of vertices in a
$k$-graph $\Lambda$ we define a $(\Lambda,X)$-family of partial
isometries subject to a set of relations similar to the
Cuntz-Krieger relations of \cite{RSY2} that generates a universal
\CSalgebra $\KGalgCLX$. We then prove a Gauge Invariant Uniqueness
Theorem (\thmref{thm:GUIT}) for $\KGalgCLX$ which generalises
\cite[Theorem 4.2]{RSY2} and then use this theorem to show that
$\KGalgCLX$ is isomorphic to a corner of $\KGalg$.

In Section 4 we describe saturated and hereditary sets of vertices
and use them to find conditions for fullness of our corners
(\corref{cor:full}). We also describe the Morita equivalence class
(\propref{prop:EquivCorner}) of a corner based on saturated and
hereditary sets. In Section 5 we examine a class of $k$-graph
morphisms which induce maps between corners and in particular how
these can be used to show Morita equivalence of different
$k$-graph algebras (\corref{cor:morph}).

In section 6 we establish necessary conditions for a corner to be
simple (\propref{prop:simple}) and in particular, if our $k$-graph
is row finite the condition is also sufficient
(\propref{prop:relaper}). In Section 7 we briefly look at some
corners that are generated by more general projections using the
dual graph defined in \cite{APS}.

Finally in Section 8 we look at skew product graphs and establish
a connection between certain fixed point algebras of $k$-graphs
and corners of skew product graphs. In particular we give a
condition for the AF core of a $k$-graph algebra to be Morita
equivalent to a skew product graph naturally associated to it.

\section{Preliminaries}

Throughout this paper we let $\N:=\set{0,1,2,\dots}$ be the set of
counting numbers and regard $\N^k$ as an abelian monoid with
identity $0=(0,0,\dots,0)$ and canonical generators
$e_i=(0,\dots,1,\dots,0)$, (1 is the $i$th coordinate). For
$n\in\N^k$ we write $n_i$ as the $i$th coordinate of $n$. There is
a partial order $\le$ on $\N^k$ given by $m\le n$ if $m_i\le n_i$
for all $1\le i\le k$, with $m<n$ if $m\le n$ and $m\ne n$. For
$m,n\in\N^k$ we write $m\vee n$ and $m\wedge n$ for their
coordinate-wise maximum and minimum respectively.

\begin{defn}
A $k$-graph is a pair $(\Lambda,d)$ consisting of a countable
category $\Lambda$ and a \emph{degree} functor
$\map{d}{\Lambda}{\N^k}$ which satisfies the \emph{factorisation
property}: for every $\lambda\in\Lambda$ and $m,n\in\N^k$ with
$d(\lambda)=m+n$ there exist unique $\mu,\nu\in\Lambda$ such that
$d(\mu)=m$, $d(\nu)=n$ and $\lambda=\mu\nu$ (see \cite{KP2} for
more details). A $k$-graph morphism is a functor between two
$k$-graphs which respects the degree map.
\end{defn}

Throughout this paper we will simply write $\Lambda$ instead of
$(\Lambda,d)$ whenever it is clear what we mean. Since we regard
$k$-graphs as a type of directed graph, we will sometimes refer to
morphisms as paths (denoted with Greek letters
$\lambda,\mu,\nu,\dots$) and objects as vertices (denote
$u,v,w,\dots$), and we will write $s$ and $r$ for the domain and
codomain maps respectively.

\begin{defn}
For all $n\in\N^k$ we define
$\Lambda^n:=\set{\lambda\in\Lambda:d(\lambda)=n}$. The
factorisation property ensures that Obj$(\Lambda)$ can be
identified with $\Lambda^0$ and we will regard them as the same
thing. Given any $v\in\Lambda^0$ and $n\in\N^k$ we define
$v\Lambda^n:=\set{\lambda\in\Lambda^n:r(\lambda)=v}$ and
$\Lambda^nv:=\set{\lambda\in\Lambda^n:s(\lambda)=v}$. Similarly,
for any $\XsubLam$, we define $X\Lambda^n:=\bigcup_{v\in
X}v\Lambda^n$ and $\Lambda^nX:=\bigcup_{v\in X}\Lambda^nv$ and
$X\Lambda:=\set{\lambda\in\Lambda:r(\lambda)\in X}$.
\end{defn}

\begin{defn}
A $k$-graph is \emph{row finite} if the set $v\Lambda^n$ is finite
for all $v\in\Lambda^0$ and $n\in\N^k$. We call a vertex
$v\in\Lambda^0$ a \emph{source} if $v\Lambda^{e_i}=\emptyset$ for
some $1\le i\le k$ and a \emph{sink} if $\Lambda^{e_i}v=\emptyset$
for some $1\le i\le k$.
\end{defn}

\begin{defn}
A $k$-graph is \emph{locally convex} if, for all $v\in\Lambda^0$
and $i,j\in\set{1,\dots,k}$ such that $i\ne j$ and
$v\Lambda^{e_i}$ and $v\Lambda^{e_j}$ are nonempty, then for all
$\lambda\in v\Lambda^{e_i}$ the set $s(\lambda)\Lambda^{e_j}$ is
nonempty.
\end{defn}


\begin{defn}
For $\lambda,\mu\in\Lambda$, we write
\begin{equation*}
\Lambdamin{\lambda}{\mu}:=\set{(\alpha,\beta):\lambda\alpha=\mu\beta,
d(\lambda\alpha)=d(\lambda)\vee d(\mu)}
\end{equation*}
for the collection of pairs which give \emph{minimal common
extensions} of $\lambda$ and $\mu$. We say the $\Lambda$ is
\emph{finitely aligned} if $\Lambdamin{\lambda}{\mu}$ is finite
(possibly empty) for all $\lambda,\mu\in\Lambda$. We also define
the (non minimal) \emph{common extensions} of $\lambda$ and $\mu$
of degree $n$ to be the set
$\Lambda^n(\lambda,\mu):=\set{(\alpha,\beta):\lambda\alpha=\mu\beta,
d(\lambda\alpha)=n}$ and note that in a finitely aligned $k$-graph
that $\Lambda^n(\lambda,\mu)$ is finite for each $n\in\N^k$ and
$(\lambda,\mu)\in\Lambda\times\Lambda$.
\end{defn}

\begin{defn} A set $E\subset v\Lambda$ is \emph{exhaustive}
if for every $\mu\in v\Lambda$ there exists $\lambda\in E$ such
that $\Lambdamin{\lambda}{\mu}\ne\emptyset$.
\end{defn}

For this paper we are only concerned with finite exhaustive sets.
This is reflected in \defref{def:kfamily}(iv) and
\defref{def:ckLXfam}(iii). We note that if $\Lambda$ is row finite
with no sources then for any $v\in\Lambda^0$ and $n\in\N^k$ the
set $v\Lambda^n$ is a finite exhaustive set.

\begin{defn}
For $n\in\N^k$, we define
\begin{equation*}
\Lambda^{\le n}:=\set{\lambda\in\Lambda:d(\lambda)\le n\text{ and
}d(\lambda)_i<n_i\Rightarrow s(\lambda)\Lambda^{e_i}=\emptyset}.
\end{equation*}
\end{defn}

If $\Lambda$ is row finite then for any $v\in\Lambda^0$ and
$n\in\N^k$ then $v\Lambda^{\le n}$ is finite exhaustive.

\begin{defn}\label{def:kfamily}
Let $(\Lambda,d)$ be a finitely aligned $k$-graph. A
\emph{Cuntz-Krieger $\Lambda$-family} is a collection
$\set{t_\lambda:\lambda\in\Lambda}$ of partial isometries in a
\CSalgebra satisfying
\begin{enumerate} \renewcommand{\theenumi}{\roman{enumi}}
\item $\set{t_v:v\in\Lambda^0}$ is a collection of mutually
orthogonal projections;

\item $t_\lambda t_\mu=t_{\lambda\mu}$ whenever $s(\lambda)=r(\mu)$;

\item $t^*_\lambda
t_\mu=\sum_{(\alpha,\beta)\in\Lambdamin{\lambda}{\mu}}t_\alpha
t^*_\beta$ for all $\lambda,\mu\in\Lambda$; and

\item $\prod_{\lambda\in E}(t_v-t_\lambda t^*_\lambda)=0$ for all
$v\in\Lambda^0$ and finite exhaustive $E\subset v\Lambda$.
\end{enumerate}
\end{defn}

\begin{rems}
\item Relation (iii) implies that $t^*_\lambda
t_\lambda=t_{s(\lambda)}$ and that $t^*_\lambda t_\mu=0$ if
$\Lambdamin{\lambda}{\mu}=\emptyset$. Also, the finitely aligned
condition is necessary for relation (iii) to make sense. See
\cite[Definition 2.5]{RSY2} for more details.

\item If $\Lambda$ is a row finite $k$-graph then relation (iv) of
\defref{def:kfamily} can be expressed as follows: for any
$v\in\Lambda^0$ and $n\in\N^k$
\begin{equation*}
t_v=\sum_{\lambda\in v\Lambda^{\le n}}t_\lambda t^*_\lambda.
\end{equation*}
\end{rems}

Given a finitely aligned $k$-graph $(\Lambda,d)$, there exists a
\CSalgebra $\KGalg$ generated by a Cuntz-Krieger $\Lambda$-family
$\set{s_\lambda:\lambda\in\Lambda}$ which is universal in the
following sense: given a Cuntz-Krieger $\Lambda$-family
$\set{t_\lambda:\lambda\in\Lambda}$ of bounded operators on a
Hilbert space $\mc{H}$, there exists a unique homomorphism
$\map{\pi}{\KGalg}{\mc{B}(\mc{H})}$ such that
$\pi(s_\lambda)=t_\lambda$ for all $\lambda\in\Lambda$. As a
consequence of \defref{def:kfamily} (i)-(iii), given
$\set{t_\lambda:\lambda\in\Lambda}$ a Cuntz-Krieger
$\Lambda$-family, by the same argument as \cite[Lemma 2.7]{RSY2},
we have
\begin{equation*}
    \KGalg=\cspan\set{t_\lambda t^*_\mu:\lambda,\mu\in\Lambda,
    s(\lambda)=s(\mu)}
\end{equation*}

Given any finitely aligned $k$-graph $(\Lambda,d)$ then it has a
strongly continuous \emph{gauge action}
$\map{\gamma}{\T^k}{\Aut(\KGalg)}$ determined by
$\gamma_z(s_\lambda)=z^{d(\lambda)}s_\lambda$ where $z\in\T^k$ and
for any $m\in\N^k$ we have $z^m=z^{m_1}_1\dots z^{m_k}_k$. The
fixed-point algebra $\KGalg^\gamma$ is AF and is equal to
$\cspan\set{s_\lambda s_\mu^*: d(\lambda)=d(\mu)}$ and is called
the \emph{AF core} of $\KGalg$ (see \cite[Theorem 3.1]{RSY2}).

\section{Cuntz-Krieger $(\Lambda,X)$-Families}

We now wish to describe corners of finitely aligned $k$-graph
algebras generated by vertex projections as a universal \CSalgebra
generated by partial isometries. A similar method was used in
\cite[\S 2]{Szy02} to describe corners of certain directed graphs.
The notation used in this paper is also comparable to that of
\cite[\S 3]{RS3} where they define rank 2 Cuntz-Krieger algebras
in a similar way.

\begin{defn}\label{def:ckLXfam}
Let $(\Lambda,d)$ be a finitely aligned $k$-graph and $\XsubLam$
be non-empty. A \emph{\CKLXfam} is a collection of partial
isometries
\begin{equation*}
    \set{T_{\alpha,\beta}:\alpha,\beta\in X\Lambda\text{ and }
    s(\alpha)=s(\beta)}
\end{equation*}
(with notation $T_\lambda:=T_{\lambda,\lambda}$ for each
$\lambda\in\Lambda$) subject to the following relations:

For any $\alpha,\beta,\lambda,\mu\in X\Lambda$ with
$s(\alpha)=s(\beta)$ and $s(\lambda)=s(\mu)$,
\begin{enumerate} \renewcommand{\theenumi}{\roman{enumi}}
\item $T_{\alpha,\beta}^*=T_{\beta,\alpha}$;
\item $T_{\alpha,\beta}T_{\lambda,\mu}=
\sum_{(\beta',\lambda')\in\Lambdamin{\beta}{\lambda}}T_{\alpha\beta',\mu\lambda'}$;
and
\item $\prod_{\lambda\in E}(T_v-T_\lambda)=0$ for all
$v\in X$ and finite exhaustive $E\subset v\Lambda$.
\end{enumerate}
\end{defn}

\begin{rems}\label{rem:defLamXfam}
\item As a result of relations (i) and (ii)
$\set{T_\lambda:\lambda\in X\Lambda}$ is a set of commuting
projections, and in particular $\set{T_v:v\in X}$ is a set of
mutually orthogonal projections.

\item When $X=\Lambda^0$ these relations reduce to a
Cuntz-Krieger $\Lambda$-family. That is, the set
$\set{T_{\lambda,s(\lambda)}:\lambda\in\Lambda}$ satisfies
\defref{def:kfamily}.
\end{rems}

 In \cite[Proposition 2.12]{RSY2} the boundary path
representation was defined for any finitely aligned $k$-graph
$\Lambda$ with Cuntz-Krieger $\Lambda$-family $\set{s_\lambda}$.
We use this representation for any \CSalgebra $B$ generated by a
\CKLXfam $\set{\Talbe}$ by noting that there is a homomorphism
$\map{\pi}{B}{\KGalg}$ given by $\pi(T_{\alpha,\beta})=s_\alpha
s^*_\beta$. Hence $B$ has a sub-representation on the boundary
path representation of $\KGalg$ and there exists a non-degenerate
\CKLXfam.

Thus we define $\KGalgCLX$ to be the universal \CSalgebra
generated by a \CKLXfam $\set{T_{\alpha,\beta}:\alpha,\beta\in
X\Lambda, s(\alpha)=s(\beta)}$. We also note, by the same argument
as \cite[Lemma 2.7(iv)]{RSY2}, that
\begin{equation*}
\KGalgCLX=\cspan\set{T_{\alpha,\beta}:\alpha,\beta\in X\Lambda,
s(\alpha)=s(\beta)}.
\end{equation*}

 By the universality of $\KGalgCLX$ and a standard
$\frac{\epsilon}{3}$ argument, $\KGalgCLX$ has a strongly
continuous \emph{gauge action} ($\gamma$) of $\T^k$ for each
$z\in\T^k$ given by,
\begin{equation*}
\gamma_z(T_{\alpha,\beta})=z^{d(\alpha)-d(\beta)}T_{\alpha,\beta}.
\end{equation*}
We call the fixed point algebra $\KGalgCLX^\gamma$ the \emph{core}
of $\KGalgCLX$.

Using a standard argument, it can be shown that
\begin{equation*}
\KGalgCLX^\gamma=\cspan\set{T_{\alpha,\beta}\in\KGalgCLX:d(\alpha)=d(\beta)}
\end{equation*}
(see \cite[lemma 3.1]{KP2}, \cite[Lemma 2.2]{BPRS} for example).

\begin{lem} \label{lem:finite}
Given a finitely aligned $k$-graph $\Lambda$ with $\XsubLam$, and
a finite set
$F=\set{\Tab{\alpha_i}{\beta_i}\in\KGalgCLX^\gamma}_{i=1}^n$. Then
there exists a finite set $\ol{F}\subset\KGalgCLX^\gamma$
containing $F$ such that $C^*(F)=\cspan\set{\ol{F}}$.
\end{lem}

\begin{proof}
Let $\Lambda(F):=\set{\alpha_i,\beta_i}_{i=1}^n$ be the set of
paths in $F$ and define $\Lambda(\le F):=
\set{\lambda\in\Lambda:\lambda\lambda'\in\Lambda(F)}$ to be the
set of initial subpaths of $F$ and let
\begin{equation*}
\Lambda(N,F):=\set{\lambda\mu\in\Lambda:\lambda\in\Lambda(\le
F),\mu\in\Lambda^n(\lambda,\nu)\text{ for some }\nu\in\Lambda(\le
F),n\le N},
\end{equation*}
where $N=\bigvee_{\lambda\in\Lambda(F)}d(\lambda)$. Then
$\Lambda(N,F)$ is the common extensions of the initial subpaths of
$F$ of degree less then or equal to $N$. Since $\Lambda$ is
finitely aligned, $\Lambda(N,F)$ is finite.

Next we note that for any
$\Talbe,\Tab{\lambda}{\mu}\in\KGalgCLX^\gamma$ then by
\defref{def:ckLXfam}(ii) their product will be the sum of elements
$\Tab{\alpha\alpha'}{\mu\mu'}\in\KGalgCLX^\gamma$ with
$d(\alpha\alpha')=d(\alpha)\vee d(\mu)$. Hence given any
$\alpha,\beta\in\Lambda(N,F)$ then for any
$(\mu,\nu)\in\Lambdamin{\alpha}{\beta}$ we must have
$\alpha\mu,\beta\nu$ in $\Lambda(N,F)$. If we let
$\ol{F}=\set{\Talbe\in\KGalgCLX^\gamma:\alpha,\beta\in\Lambda(N,F)}$
then $\ol{F}$ contains $F$, is closed under adjoints and
$\cspan\set{\ol{F}}$ is closed under multiplication.
\end{proof}

\begin{lem}\label{lem:Afcore}
Let $\Lambda$ be a finitely aligned $k$-graph with $\XsubLam$.
Then the fixed point algebra $\KGalgCLX^\gamma$ is AF.
\end{lem}

\begin{proof}
By \cite[Theorem 2.2]{Br} it suffices to show that for any finite
set $F\subset \KGalgCLX^\gamma$, that $\CS{F}$ is finite
dimensional. Without loss of generality, we may assume that
$F=\set{\Tab{\alpha_i}{\beta_i}\in\KGalgCLX^\gamma}_{i=1}^n$ and
by \lemref{lem:finite} there exists a finite set $\ol{F}$ such
that span$\set{\ol{F}}$ is closed under multiplication and taking
adjoints and such that $\CS{F}\subset\cspan\set{\ol{F}}$.
\end{proof}

 We now wish to prove an analogue of the Gauge Invariant
Uniqueness Theorem for \CSalgebras generated by \CKLXfams[].

\begin{thm}[Gauge Invariant Uniqueness Theorem]\label{thm:GUIT}
Let $\Lambda$ be a finitely aligned $k$-graph with $\XsubLam$,
$\set{t_{\alpha,\beta}}$ a \CKLXfam and $\pi$ be a representation
of $\KGalgCLX$ such that $\pi(\Talbe)=t_{\alpha,\beta}$. Suppose
that for each $v\in \Lambda^0$ with $X\Lambda v\ne\emptyset$ there
exists a path $\chi_v\in X\Lambda v$ such that
$\pi(\Tl{\chi_v})\ne0$, and suppose that there is a strongly
continuous action $\delta$ of $\T^k$ on $\CS{t_{\alpha,\beta}}$
such that $\delta_z\circ\pi=\pi\circ\gamma_z$ for all $z\in\T^k$.
Then $\pi$ is faithful.
\end{thm}

\begin{proof}

First suppose  $\Talbe\in\KGalgCLX^\gamma$ and $\pi(\Talbe)=0$,
then since
\begin{equation*}
\Tl{\chi_{s(\alpha)}}=\Tab{\chi_{s(\alpha)}}{\alpha}\left(
\Talbe\right)\Tab{\beta}{\chi_{s(\alpha)}}\ ,\tag{$*$}\label{eq:*}
\end{equation*}
we must have that $\pi(\Tl{\chi_{s(\alpha)}})=0$ which contradicts
our hypothesis.

Next suppose $x=\sum_{i=1}^nc_i\Tab{\alpha_i}{\beta_i}\in
\KGalgCLX^\gamma$ with $\pi(x)=0$ and further suppose that
$\Lambda$ has no sources. Let $N=\bigvee_{i=1}^n d(\alpha_i)$.
There exists $\lambda,\mu\in\Lambda^N$ such that
$\lambda=\alpha_j\lambda'$ and $\mu=\beta_j\lambda'$ for some
$1\le j\le n$. Hence $d(\lambda)\ge d(\alpha_i)$ and $d(\mu)\ge
d(\beta_i)$ (since $d(\alpha_i)=d(\beta_i)$) for all $1\le i\le
n$. Thus if $\Lambdamin{\lambda}{\alpha_i}\ne\emptyset$ then
$\lambda=\alpha_i\lambda''$ and similarly $\mu=\beta_i\mu''$
whenever $\Lambdamin{\mu}{\beta_i}\ne\emptyset$.Hence
\begin{equation*}
\Tlam\left(c_i\Tab{\alpha_i}{\beta_i}\right)\Tab{\mu}{\lambda}=
\begin{cases}c_i\Tlam&i=j\OR\Lambdamin{\lambda}{\alpha_i}\ne\emptyset,\ \Lambdamin{\mu}{\beta_i}\ne\emptyset\\
0&\\\end{cases}
\end{equation*}
and so $\Tlam.x.\Tab{\mu}{\lambda}=c_x\Tlam\ne0$. Hence $\pi(x)=0$
implies $\pi(\Tlam)=0$ which implies
$\pi(\Tl{\chi_{s(\lambda)}})=0$ by (\ref{eq:*}). Hence $\pi(x)\ne
0$.

Next suppose that $\Lambda$ has sources. There still exists
$\lambda,\mu\in\Lambda^{\le N}$ such that
$\lambda=\alpha_j\lambda'$ and $\mu=\beta_j\lambda'$ for some
$1\le j\le n$ with the property that for all $1\le i\le n$ then
$d(\lambda)\ngeq d(\alpha_i)$ implies
$\Lambdamin{\lambda}{\alpha_i}=\emptyset$ and (since
$d(\alpha_i)=d(\beta_i)$) the same for $\mu$ with each $\beta_i$.
Hence we still have $\Tlam.x.\Tab{\mu}{\lambda}=c_x\Tlam$ as
before.

Hence $\pi$ is faithful on
$\mc{F}=\spn\set{\Talbe\in\KGalgCLX:d(\alpha)=d(\beta)}$ and since
$\KGalgCLX^\gamma$ is AF by \lemref{lem:Afcore} then every non
trivial ideal in $\KGalgCLX^\gamma$ must intersect $\mc{F}$ by
\cite[Lemma 3.1]{Br} and since the kernel of $\pi$ is an ideal
then $\pi$ must also be faithful on $\KGalgCLX^\gamma$. Finally,
since $\pi$ is faithful on $\KGalgCLX^\gamma$ which is AF, the
remainder of the proof is now standard (see \cite[Theorem
3.4]{KP2} or \cite[Theorem 4.2]{RSY2})
\end{proof}

\begin{rem}
If $X=\Lambda^0$ then for each $v\in\Lambda^0$ we may take
$\chi_v=v$, and then \thmref{thm:GUIT} becomes the usual Gauge
Invariant Uniqueness Theorem for finitely aligned $k$-graphs as
seem in \cite[Theorem 4.2]{RSY2}.
\end{rem}

Given a $k$-graph $\Lambda$ with $\XsubLam$ and Cuntz-Krieger
$\Lambda$-family $\set{s_\lambda}$, then by the same argument as
\cite[Lemma 3.3.1]{PR} the sum $\sum_{v\in X}s_v$ converges to a
projection $P_X\in\mc{M}(\KGalg)$.

\begin{cor}\label{cor:Corner}
Let $\Lambda$ be a finitely aligned $k$-graph and
$\set{s_\lambda}$ be a Cuntz-Krieger $\Lambda$-family. Let
$\XsubLam$ and $\set{\Talbe}$ be a \CKLXfam then
\begin{equation*}
\Corner{\Lambda}{X}\cong \KGalgCLX.
\end{equation*}
\end{cor}

\begin{proof} Define a map $\map{\phi}{\KGalgCLX}{\Corner{\Lambda}{X}}$
by $\phi(\Talbe)=s_\alpha s^*_\beta$. Then $\phi$ is a surjective
homomorphism such that $\gamma_z(\phi(\Talbe))=\gamma_z s_\alpha
s^*_\beta=\phi(\gamma_z(T_{\alpha,\beta}))$. Since
$\phi(\Tlam)=s_\lambda s^*_\lambda\ne0$ for all $\lambda\in
X\Lambda v$ then by \thmref{thm:GUIT}, $\phi$ is also injective.
\end{proof}

\begin{rem}\label{rem:thesame}
Using the map $\phi$ in the proof of \corref{cor:Corner} we have
that
\begin{equation*}
\KGalgC{\Lambda^0}\cong \KGalg.
\end{equation*}
In this case the relations in \defref{def:ckLXfam} are equivalent
to the relations of a Cuntz-Krieger $\Lambda$-family as given in
\defref{def:kfamily}. Hence when it is convenient we
will identify $\KGalgC{\Lambda^0}$ with $\KGalg$ via the mapping
$\Talbe\mapsto s_\alpha s^*_\beta$ in order to avoid conflicts of
notation.
\end{rem}


\section{Fullness of $\KGalgCLX$}

When considering corners it is natural to ask when the corner is
full. The answer has a lot to do with saturated hereditary subsets
of $\Lambda^0$ and their association with gauge invariant ideals
in $\KGalgCLL$ (see \cite{Si2} for details).

\begin{defn}
Given a $k$-graph $\Lambda$ with $H,S\subset\Lambda^0$ then:
\begin{enumerate}
\item we say $H$ is \emph{hereditary} if for all $v\in H$ then
$v\Lambda=v\Lambda H$,
\item we say $S$ is \emph{saturated} if for any $v\in \Lambda^0$
such that there exists a finite exhaustive set $E\subset v\Lambda
S$ then $v\in S$.
\end{enumerate}
If $H$ is the smallest hereditary set containing $V$ and $S$ is
the smallest saturated set containing $H$, then $S$ is the
smallest saturated and hereditary set containing $V$ (see
\cite[Lemma 3.2]{Si2}). We call the smallest saturated hereditary
set containing $V$ the \emph{saturation} of $V$ and denote it as
$\sat{V}$.
\end{defn}

\begin{prop}\label{prop:EquivCorner}
Given a finitely aligned $k$-graph $\Lambda$ with $\XsubLam$ then
$\KGalgCLX$ is Morita equivalent to $\KGalgC{\sat{X}}$.
\end{prop}

\begin{proof}
By \cite[Example 3.6]{RW} $\KGalgCLX$ is Morita equivalent to the
ideal generated by $P_X$ and by \cite[Lemma 3.3]{Si2} this ideal
is equal to the ideal generated by $P_{\sat{X}}$.
\end{proof}

\begin{rems}
\item If $\Lambda$ is finitely aligned and $\XsubLam$ then
$\KGalgCLX$ is not usually a graph algebra but is always Morita
equivalent to a graph algebra because of
\propref{prop:EquivCorner}. If $\Lambda$ is row finite and locally
convex and $\XsubLam$ is a hereditary set then by \cite[Theorem
5.2(c)]{RSY1} $\KGalgCLX\cong\KGalgC[X\Lambda]{X}$ where
$X\Lambda$ is a $k$-graph because $X$ is hereditary. Since
$X\Lambda^0=X$ then it follows that $\KGalgC[X\Lambda]{X}$ is a
graph algebra. If $\Lambda$ is finitely aligned but not row finite
and $\XsubLam$ is a saturated and hereditary set then by
\cite[Lemma 3.6]{Si2} we again have
$\KGalgCLX\cong\KGalgC[X\Lambda]{X}$. Hence for any finitely
aligned $k$-graph with $\XsubLam$ then by
\propref{prop:EquivCorner} we have $\KGalgCLX$ is Morita
equivalent to $\KGalgC{\sat{X}}$ with the latter being a graph
algebra.

\item Since Morita equivalence respects many \CSalgebra
properties (e.g. simplicity, AF, etc) then in many cases where it
may be convenient we may assume without loss of any generality
that our set $\XsubLam$ is saturated and hereditary.

\end{rems}

In particular, \propref{prop:EquivCorner} says that $\KGalgCLX$ is
a full corner of $\KGalgC{\sat{X}}$, and hence we have
\corref{cor:full}

\begin{cor}\label{cor:full} Given a finitely aligned $k$-graph
$\Lambda$ with $\XsubLam$ then $\KGalgCLX$ is a full corner of
$\KGalgCLL$ if and only if $\sat{X}=\Lambda^0$.
\end{cor}

\begin{proof}
By \propref{prop:EquivCorner} if $\sat{X}=\Lambda^0$ then
$\KGalgCLX$ is full in $\KGalgCLL$. Conversely, if $\KGalgCLX$ is
full then the ideal generated by $P_X$ is equal to the ideal
generated by $P_{\Lambda^0}$ and hence $\sat{X}=\Lambda^0$.
\end{proof}

\begin{cor}
Let $\Lambda$ be a finitely aligned $k$-graph and let
$X,Y\subseteq\Lambda^0$ be such that $\sat{X}=\sat{Y}$. Then
$\KGalgCLX$ is Morita equivalent to $\KGalgC{Y}$.
\end{cor}

\begin{proof} $\KGalgCLX$ is Mortia
equivalent to $\KGalgC{\sat{X}}=\KGalgC{\sat{Y}}$ which is Morita
equivalent to $\KGalgC{Y}$ by \propref{prop:EquivCorner}.
\end{proof}

\begin{eg}
Let $\XsubLam$ be any subset and let $X^c=\Lambda^0\setminus X$
and suppose $\sat{X}=\sat{X^c}$. This implies $\sat{X}=\Lambda^0$
and hence $\KGalgCLX$ and $\KGalgC{X^c}$ are complimentary corners
(see \cite[Theorem 1.1]{BGR}).
\end{eg}

\section{$k$-graph Morphisms}
\begin{defn} Given $k$-graphs $\Lambda_1$, $\Lambda_2$ with
$\XsubLam$ and a $k$-graph morphism
$\map{\phi}{\Lambda_1}{\Lambda_2}$ then we say $\phi$ is
\emph{saturated} with respect to $X$ if
$\map{\phi}{X\Lambda_1}{\phi(X)\Lambda_2}$ is a bijection (c.f.
\cite[Definition 3.2 \& Proposition 3.3]{PRY}). If $X=\Lambda^0$
then we call $\phi$ a \emph{saturated} $k$-graph morphism.
\end{defn}

\begin{thm}\label{thm:morph}
Given a finitely aligned $k$-graphs $\Lambda_1$, $\Lambda_2$ with
$X\subseteq\Lambda_1^0$ and a $k$-graph morphism
$\map{\phi}{\Lambda_1}{\Lambda_2}$ that is relatively saturated
with respect to $X$ then
$\KGalgC[\Lambda_1]{X}\cong\KGalgC[\Lambda_2]{\phi(X)}$.
\end{thm}

\begin{proof}
This proof follows the same argument as \cite[Proposition
3.3]{PRY} and is repeated here for convenience. Let $\set{\Talbe}$
be a Cuntz-Krieger $(\Lambda_1,X)$-family and let
$\set{\Tab[S]{\gamma}{\delta}}$ be a Cuntz-Krieger
$(\Lambda_2,\phi(X))$-family. The relative saturation property
ensures that $\set{\Tab[S]{\phi(\alpha)}{\phi(\beta)}}$ is a
Cuntz-Krieger $(\Lambda_2,\phi(X))$-family and the universal
property of $\KGalgC[\Lambda_1]{X}$ induces a homomorphism
$\map{\phi_*}{\KGalgC[\Lambda_1]{X}}{\KGalgC[\Lambda_2]{\phi(X)}}$
given by $\phi_*(\Talbe)=\Tab[S]{\phi(\alpha)}{\phi(\beta)}$. Then
$\phi_*$ is surjective since $\phi$ is saturated and also $\phi$
is degree preserving since it is a $k$-graph morphism and hence
$\phi_*$ respects the gauge action. Finally
$\phi_*(\Tlam)=\Tl[S]{\phi(\lambda)}\ne0$ so by \thmref{thm:GUIT}
$\phi_*$ is also injective.
\end{proof}

\begin{eg}\label{eg:hereditary}
Let $\Lambda$ be a finitely aligned $k$-graph with
$H\subset\Lambda$ a hereditary subset. Then $H\Lambda$ is a sub
$k$-graph of $\Lambda$ and the identity map
$\map{i}{H\Lambda}{\Lambda}$ is a relatively saturated $k$-graph
morphism with respect to $H$. Hence we have
$\KGalg[H\Lambda]\cong\KGalgC{H}$ by \thmref{thm:morph} which is
an improvement of \cite[Theorem 5.2(c)]{RSY1} which requires
$\Lambda$ to be row finite locally convex and \cite[Lemma
3.6]{Si2} which requires $H$ to be saturated and hereditary.
\end{eg}

\begin{cor}\label{cor:morph}
Given finitely aligned $k$-graphs $\Lambda_1$, $\Lambda_2$ and a
saturated $k$-graph morphism $\map{\phi}{\Lambda_1}{\Lambda_2}$
then $\KGalg[\Lambda_1]$ is Morita equivalent to
$\KGalg[\Lambda_2]$ if and only if
$\sat{\phi(\Lambda_1^0)}=\Lambda_2^0$.
\end{cor}

\begin{proof}
Follows from \thmref{thm:morph} and \corref{cor:full}.
\end{proof}

\begin{eg}\label{eg:OmegaDelta} Recall from \cite{KP3} that
$\Omega_k$ is the $k$-graph with vertex set $\N^k$ and paths
$\set{(m,n)\in\N^k\times\N^k:m<n}$ with $r(m,n)=m$, $s(m,n)=n$ and
$d(m,n)=m-n$, while $\Delta_k$ is the $k$-graph with vertices
$\Z^k$ and paths $\set{(m,n)\in\Z^k\times\Z^k:m<n}$ with the same
range, source and degree maps as above. So there is a natural
embedding using the identity map $\map{i}{\Omega_k}{\Delta_k}$
which is a saturated $k$-graph morphism. Noting that
$\Omega_k^0=\N^k$ is a hereditary subset of $\Delta_k^0=\Z^k$
then, as in \egref{eg:hereditary},
$\KGalg[\Omega_k]\cong\KGalgC[\Delta_k]{\N^k}$. However
$\sat{\N^k}=\Z^k$ so by \corref{cor:morph}
$\KGalgC[\Delta_k]{\N^k}$ is a full corner of $\KGalg[\Delta_k]$
and thus we have $\KGalg[\Omega_k]$ is Morita equivalent to
$\KGalg[\Delta_k]$.
\end{eg}

We now look at an application of saturated $k$-graph morphisms
with symbolic dynamics in which we are concerned with the
bi-infinite path space of a $k$-graph.

\begin{defn}
Given a $k$-graph $\Lambda$ the \emph{bi-infinite path space} of
$\Lambda$ is the set of $k$-graph morphisms
\begin{equation*}
\Lambda^\Delta=\set{\map{x}{\Delta_k}{\Lambda}}
\end{equation*}
where $\Delta_k$ is as defined in example \ref{eg:OmegaDelta}.
\end{defn}

A typical construction is the \emph{essential subgraph} which is
the largest subgraph with no sinks or sources. The bi-infinite
path space of the essential subgraph can be identified with the
bi-infinite path space of the original graph (see \cite{LM}) and
is also identified with the edge shift associated to the graph.
This was done for 1-graphs in \cite{LM}.

We will now adapt this construction to $k$-graphs and show
conditions for a $k$-graph and its essential subgraph to have
Morita equivalent algebras. Constructing the essential subgraph
involves removing `stranded' vertices that do not lie on any
bi-infinite path.

\begin{rem} If $\Lambda^\Delta=\emptyset$ then the essential
subgraph will be trivial.
\end{rem}

\begin{defn}
Given a finitely aligned $k$-graph $\Lambda$ then for any
$v\in\Lambda^0$ we say $v$ is \emph{stranded} if there exists
$n\in\N^k$ such that $v\Lambda^n=\emptyset$ or
$\Lambda^nv=\emptyset$. We denote $\StL$ as the set of all
stranded vertices in $\Lambda^0$.
\end{defn}

For any $k$-graph $\Lambda$ we construct the essential subgraph
$\eLambda$ by first identifying all the stranded vertices in
$\Lambda^0$ and then constructing the subcategory with objects
$\Lambda^0\setminus\StL$ and morphisms
$\set{\lambda\in\Lambda:s(\lambda),r(\lambda)\notin\StL}$ (cf.
\cite{LM}). We can constructively define the set of stranded
vertices in a recursive manner (c.f. \cite[Remark 3.1]{BHRS}) as
follows. First let $S_0$ be the set of all sinks and sources. Next
let
\begin{equation*}
S_{n+1}=S_n\bigcup_{i=1}^k\set{v\in\Lambda^0:v\Lambda^{e_i}\subset
v\Lambda S_n}\cup\set{v\in\Lambda^0:\Lambda^{e_i}v\subset
S_n\Lambda v},
\end{equation*}
and finally let $\StL=\bigcup_{n\in\N}S_n$.

It is worth taking a moment to check that $\eLambda$ forms a valid
$k$-graph.

\begin{lem}
Let $\Lambda$ be a finitely aligned $k$-graph with
$\Lambda^\Delta\ne\emptyset$ then there exists a unique
non-trivial $k$-graph $\eLambda\subset\Lambda$ such that
$\eLambda$ the largest subgraph of $\Lambda$ with no sinks or
sources and $\eLambda^\Delta=\Lambda^\Delta$.
\end{lem}

\begin{proof}
Let $\eLambda$ be the subcategory of $\Lambda$ with objects
$\Lambda^0\setminus\StL$ and morphisms
$\set{\lambda\in\Lambda:s(\lambda),r(\lambda)\in\eLambda^0}$.
Clearly $\eLambda$ is nontrivial if $\Lambda^\Delta\ne\emptyset$.
To see that $\eLambda$ is a $k$-graph we need to check the
factorisation property. Suppose $\lambda\in\eLambda$, then this
implies that $\lambda\in\Lambda$ and $r(\lambda)\Lambda^n$ and
$\Lambda^ns(\lambda)$ are nonempty for all $n\in\N^k$. Let
$p,q\in\N^k$ such that $p+q=d(\lambda)$ then there exist
$\mu\in\Lambda^p$ and $\nu\in\Lambda^q$ such that
$\lambda=\mu\nu$. By the factorisation property of $\Lambda$ we
must have $\Lambda^ns(\mu)$ and $s(\mu)\Lambda^m$ nonempty for all
$m\le p$ and $n\le q$ and further since $s(\lambda)$ and
$r(\lambda)$ are not stranded we then have $s(\mu)\Lambda^n$ and
$\Lambda^ns(\mu)$ are nonempty for all $n\in\N^k$ and hence
$\mu,\nu\in\eLambda$. Therefore $\eLambda$ is a $k$-graph.

Clearly if $x\in\Lambda^\Delta$ then $x\in\eLambda^\Delta$ and
vice-versa so $\Lambda^\Delta=\eLambda^\Delta$ and clearly
$\eLambda$ is unique. To show $\eLambda$ is the largest subgraph
with no sinks or sources, suppose $S$ is also s subgraph of
$\Lambda$ with no sinks or sources. Then every $v\in S^0$ is not
stranded and hence $S\subset\eLambda$.
\end{proof}

\begin{eg}
Let $\Omega_k$ and $\Delta_k$ be defined as in
\egref{eg:OmegaDelta}. Then $\ess{\Omega_k}=\emptyset$ since every
vertex is stranded, however $\ess{\Delta_k}=\Delta_k$ since every
vertex is not stranded.
\end{eg}

\begin{defn}
A $k$-graph $\Lambda$ is \emph{essentially saturated} if for every
$v\in\Lambda^0$ there exists $x\in\Lambda^\Delta$ and $n\in\Z^k$
such that $v\Lambda x(n)\ne\emptyset$.
\end{defn}

The essentially saturated property says that for every vertex
$v\in\Lambda^0$ there exists a path $\lambda\in
v\Lambda\eLambda^0$ between $v$ and the the essential subgraph
$\eLambda$. Further more it implies that $\Lambda$ has no sources
and $\Lambda^\Delta\ne\emptyset$ and in particular for every
$v\in\Lambda^0$ then $v\Lambda^n$ is non empty for all $n\in\N^k$
and that there exists $w\in\Lambda^0$ such that $v\Lambda
w\ne\emptyset$ and $\Lambda^nw$ is nonempty for all $n\in\N^k$
(possibly with $v=w$).

\begin{defn}\label{def:FE} A $k$-graph $\Lambda$ is \emph{finitely
exhaustive} if for every $v\in\Lambda^0$ there exists a finite
exhaustive set $E\subset v\Lambda$.
\end{defn}

Any row finite $k$-graph is automatically finitely exhaustive,
however there does exist a class of finitely exhaustive $k$-graphs
that are not row finite. In many ways these $k$-graph would behave
like row finite $k$-graphs because \defref{def:ckLXfam} allows
each vertex projection to be written as a finite sum of path
projections. We will consider finitely exhaustive $k$-graphs again
in section 8.

\begin{lem}\label{lem:essSat}
Given a finitely exhaustive $k$-graph $\Lambda$ that is
essentially saturated then $\KGalg[\eLambda]$ is Morita equivalent
to $\KGalg$.
\end{lem}

\begin{proof} If $\Lambda$ is essentially saturated then
$\eLambda^0$ is a hereditary subset of $\Lambda^0$ and since
$\Lambda$ is finitely exhaustive $\sat{\eLambda^0}=\Lambda^0$.
Hence the result follows from \corref{cor:morph} by the same
argument as Example \ref{eg:hereditary}.
\end{proof}

\section{Simplicity}
\begin{defn}\label{def:relCof} A $k$-graph is \emph{cofinal} if
$\sat{\set{v}}=\Lambda^0$ for all $v\in\Lambda^0$. Also, given
$\XsubLam$ then $X$ is \emph{relatively cofinal} if
$\sat{\set{v}}=\sat{X}$ for all $v\in X$.
\end{defn}

Cofinal is usually defined using the infinite path space of
$\Lambda$ (see \cite[Definition 4.7]{KP2} or \cite[Definition
8.4]{Si2}). In short a $k$-graph is cofinal if for every vertex in
$\Lambda$ and every infinite path in $\Lambda^\infty$ there exists
a finite path connecting them, while the relatively cofinal
condition says that every vertex in $X$ can be connected to every
path in $\Lambda^\infty$. It should also be noted that $\Lambda$
is cofinal if and only if $\Lambda^0$ is relatively cofinal.

\begin{prop}\label{prop:simple}
Let $\Lambda$ be a finitely aligned $k$-graph with $\XsubLam$ such
that all ideals in $\KGalgCLX$ are gauge invariant, then
$\KGalgCLX$ is simple if and only if and $X$ is relatively cofinal
(c.f. \cite[Theorem 12]{Szy01}).
\end{prop}

\begin{proof}
By \propref{prop:EquivCorner} we may assume $X$ is saturated and
hereditary and so $\KGalgCLX\cong\KGalg[X\Lambda]$ by
\thmref{thm:morph}. In particular, if $X$ is relatively cofinal in
$\Lambda$ then $X\Lambda$ is cofinal. Finally by \cite[Proposition
8.5]{Si2} $\KGalg[X\Lambda]$ is simple if and only if $X\Lambda$
is cofinal.
\end{proof}

For most purposes \propref{prop:simple} is unsatisfactory for
determining simplicity since there is not yet a necessary and
sufficient condition for the ideals of $k$-graph to all be gauge
invariant. In \cite[Theorem 7.2]{Si2} condition (D) is stated for
when all the ideals of a finitely aligned $k$-graph are gauge
invariant however it is not easily checkable. However for row
finite $k$-graphs we can say much more.\\

Recall from \cite[Definition 4.1]{KP2} that $x\in\Lambda^\infty$
is periodic if there exists $p\in\Z^k$ such that
$x(m+p,n+p)=x(m,n)$ for all $m,n\in\N^k$ (with $m+p\ge 0$) and is
eventually periodic if there exists $n\in\N^k$ such that
$\sigma^{n}x$ is periodic (where $\sigma$ is the shift map). A
path in $\Lambda^\infty$ is \emph{aperiodic} if it is not periodic
or eventually periodic.

\begin{defn}
Let $\Lambda$ be a $k$-graph with $\XsubLam$. Then we say $X$ is
\emph{relatively aperiodic} if for all $v\in X$ there exists $x\in
v\Lambda^\infty$ such that $x$ is aperiodic. We say $\Lambda$ is
\emph{aperiodic} if $\Lambda^0$ is relatively aperiodic.
\end{defn}

In \cite[Proposition 4.8]{KP2} it was shown that if $\Lambda$ is a
row finite $k$-graph then $\KGalg$ is simple if $\Lambda$ is
aperiodic and cofinal. Here we simply extend this notion to
corners of locally convex and row finite $k$-graphs by using the
definition of relative aperiodicity.

\begin{prop}\label{prop:relaper}
Let $\Lambda$ be a row finite $k$-graph with $\XsubLam$ and let
$X$ be relatively aperiodic. Then $\KGalgCLX$ is simple if and
only if $X$ is relatively cofinal.
\end{prop}

\begin{proof} By \propref{prop:EquivCorner} we may assume $X$
is saturated and hereditary and $\KGalgCLX\cong\KGalg[X\Lambda]$
by \thmref{thm:morph}. We note that if $X$ is relatively aperiodic
in $\Lambda$ then $X\Lambda$ is aperiodic. Hence \cite[Proposition
4.8]{KP2} applies and $\KGalgCLX$ is simple if and only if
$X\Lambda$ is cofinal. Since $X$ is relatively cofinal if and only
if $X\Lambda$ is cofinal then this completes the proof.
\end{proof}

\section{Corners Generated by Subsets of $\Lambda^p$}

As a consequence of \propref{prop:relaper} if $\Lambda$ is a row
finite $k$-graph and $\XsubLam$ is relatively aperiodic then all
ideals in $\KGalgCLX$ are generated by saturated hereditary
subsets of $X$. In such a case all corners of $\KGalgCLL$ are
generated by subsets of $\Lambda^0$ and we need only consider $X$
as a subset of $\Lambda^0$. However if $X$ is not relatively
aperiodic not not all ideals are generated by saturated hereditary
subsets of $\Lambda^0$ (see\cite[\S 5]{Si2}) and hence not all
corners can be generated by subsets of $\Lambda^0$. In this
section we show that if $\Lambda$ is any row finite $k$-graph with
no sources then we can easily extend our ideas in this paper to
corners generated by certain subsets of $\Lambda^p$ for some
$p\in\N^k$.

\begin{defn}
Let $\Lambda$ be a $k$-graph and let $p\in\N^k$, then the
\emph{dual graph} is
$p\Lambda:=\set{\lambda\in\Lambda:d(\lambda)\ge p}$ with range,
source and degree maps defined as follows: For any $\lambda\in
p\Lambda$ with $\lambda=\sigma\lambda'=\lambda''\rho$ and
$d(\sigma)=d(\rho)=p$,
\begin{equation*}
r_p(\lambda)=\rho,\quad s_p(\lambda)=\sigma,\quad
d_p(\lambda)=d(\lambda)-p.
\end{equation*}
and composition defined as follows: For any
$\lambda=\lambda'\rho,\mu=\rho\mu'\in p\Lambda$ with
$r_p(\lambda)=s_p(\mu)=\rho$, then
$\lambda\circ_p\mu=\lambda'\rho\mu'$.
\end{defn}

 For more details of dual higher rank graphs see \cite[\S
3]{APS}. In particular $p\Lambda$ is a $k$-graph and if $\Lambda$
is row finite with no sources then $\KGalg\cong\KGalg[p\Lambda]$
(see \cite[Proposition 3.2 \& Theorem 3.5]{APS}).

\begin{lem}\label{lem:dual}
Let $\Lambda$ be a row finite $k$-graph with no sources and let
$X\subseteq(p\Lambda)^0=\Lambda^p$ for some $p\in\N^k$ be such
that for any $\alpha,\beta\in X$ we have
$\Lambdamin{\alpha}{\beta}=\emptyset$. Then
\begin{equation*}
P_X\KGalg P_X\cong\KGalgC[p\Lambda]{X},
\end{equation*}
where $P_X=\sum_{\lambda\in X}s_\lambda s^*_\lambda$.
\end{lem}

\begin{proof}
The hypothesis on the set $X$ ensures that $\set{s_\lambda
s^*_\lambda:\lambda\in X}$ is a set of mutually orthogonal
projections in $\KGalg$ so $P_X$ converges to a projection in
$\mc{M}\left(\KGalg\right)$ (by the same argument as \cite[Lemma
3.3.1]{PR}). The rest follows from \corref{cor:Corner} and
\cite[Theorem 3.5]{APS}.
\end{proof}

We can extend \defref{def:ckLXfam} to include
$X\subseteq\Lambda^p$ subject to the hypothesis in
\lemref{lem:dual}. Note that if $p=0$ then \lemref{lem:dual}
reduces to \corref{cor:Corner}.

When we talk about corners generated by arbitrary subsets of
$\Lambda^p$ we have to be careful to watch that $P_X$ converges to
a projection in $\mc{M}\left(\KGalg\right)$. However in some cases
we can still talk about \CKLXfams generated by such arbitrary
subsets. For example let us suppose $\Lambda$ is row finite with
no sinks or sources with $X\subset\Lambda^p$ and suppose
$\mu,\nu\in X$ with $\mu=\nu\nu'$. By definition
\begin{equation*}
\KGalgC[p\Lambda]{X}=\cspan\set{\Talbe:\alpha,\beta\in
X\Lambda,s_p(\alpha)=s_p(\beta)},
\end{equation*}
and hence for any $\lambda\in\mu\Lambda\cup\nu\Lambda$ there is a
$\Tlam\in\KGalgC[p\Lambda]{X}$. However
$\mu\Lambda\subset\nu\Lambda$ so we would have
$\KGalgC[p\Lambda]{X\setminus\set{\mu}}=\KGalgC[p\Lambda]{X}$ is
the sense that they are isomorphic via the identity map.\\

\section{Skew Product Graphs}

In this section we look at a $k$-graph construction called a skew
product graph $\SPg{G}{\Lambda}$ and its relation to fixed point
algebras.

\begin{defn}
Given a $k$-graph $\Lambda$ and a functor $\map{c}{\Lambda}{G}$
where $G$ is a discrete group then the \emph{skew product} graph
$\SPgGL$ is the $k$-graph with objects $\SPg[]{G}{\Lambda^0}$ and
morphisms $\SPg[]{G}{\Lambda}$ with
$s(g,\lambda)=(gc(\lambda),s(\lambda))$ and
$r(g,\lambda)=(g,r(\lambda))$ and degree map
$d(g,\lambda)=d(\lambda)$ (see \cite[Definition 5.1]{KP2} for
details).
\end{defn}

In particular, a functor $\map{c}{\Lambda}{G}$ gives rise to a
normal coaction $\gamma_c$ of $G$ on $\KGalgCLL$ given by:
\begin{equation*}
\gamma_c(\Talbe)=\Talbe\otimes 1_{c(\alpha)c(\beta)^{-1}},
\end{equation*}
 where for any $g\in G$ then $1_g$ is the point mass
function in $C^*(G)$. Then the fixed-point algebra is
$\KGalgCLL^{\gamma_c}=\cspan\set{\Talbe\in\KGalgCLL:c(\alpha)=c(\beta)}$
(see \cite[\S 7]{PQR2} for more details).

\begin{prop}\label{prop:AFcore}
Let $\Lambda$ be a finitely aligned $k$-graph, let $G$ be a
discrete group with a functor $\map{c}{\Lambda}{G}$ and let
$\gamma_c$ be the corresponding coaction of $G$ on $\KGalgCLL$.
Then
\begin{equation*}
\KGalgCLL^{\gamma_c}\cong\KGalgCGLV,
\end{equation*}
where $V=\set{(1,v)\in G\times\Lambda^0}$ and 1 is the identity
element of $G$.
\end{prop}

\begin{proof}
For any $(g,\lambda)\in\SPgGL$ with range in $V$ we must have
$g=1$. Also for any $(1,\mu),(1,\nu)\in\SPgGL$ with the same
source we must have $c(\mu)=c(\nu)$ and $s(\mu)=s(\nu)$. Hence
$\CS{\SPgGL}=\cspan\set{\Tab{(1,\mu)}{(1,\nu)}:c(\mu)=c(\nu),s(\mu)=s(\nu)}$.
Thus we define the map $\map{\phi}{\KGalgCGLV}{\KGalgCLL^\gamma}$
by $\phi\left(\Tab{(1,\mu)}{(1,\nu)}\right)=\Tab{\mu}{\nu}$ which
is clearly a surjective homomorphism. Since
$\phi\left(\Tl{(1,\lambda)}\right)=\Tlam\ne0$ for all
$\lambda\in\Lambda$ and $\phi$ is gauge invariant then by
\thmref{thm:GUIT} $\phi$ is also injective.
\end{proof}

 \propref{prop:AFcore} is a well known property of directed
graphs (see \cite[Theorem 4.6]{Cr} and \cite[Proposition 2.8]{KP1}
to name a few) and as such, it is no surprise that it is also true
for finitely-aligned $k$-graphs. However the main difference here
is by using the universal property of corner algebras and
\thmref{thm:GUIT} we obtain a short proof.

\begin{eg} For any finitely aligned $k$-graph $\Lambda$, let
$G=\Z^k$ and take the degree map as our functor. Then
$\KGalgCZLV\cong\KGalgCLL^{\gamma_d}$ by \propref{prop:AFcore}
where $\gamma_d=\gamma$ is the gauge action.
\end{eg}

In particular, if $\KGalgCZLV$ is a full corner of $\CS{\SPgZL}$,
then the $\KGalgCLL^\gamma$ is Morita equivalent to the
$\CS{\SPgZL}$. So our next aim is to find a condition for the
fullness of $\KGalgCZLV$.

\begin{prop}\label{prop:MEcore}
Let $\Lambda$ be a finitely aligned $k$-graph and
$V=\set{(0,v)\in\Z\times\Lambda^0}$ then $\KGalgCZLV$ is a full
corner of $\CS{\SPgZL}$ if and only if $\Lambda$ is essentially
saturated and finitely exhaustive.
\end{prop}

\begin{proof}
Suppose $\Lambda$ is essentially saturated and finitely
exhaustive. Then for any $(x,v),(y,w)\in\SPgZL^0$ then there
exists a path in $(x,v)\left(\SPgZL\right)(y,w)$ if and only if
there exists a path $\lambda\in v\Lambda w$. Thus $\SPgZL$ is
essentially saturated if and only if $\Lambda$ is essentially
saturated, so by using \lemref{lem:essSat} we may also assume
$\SPgZL$ has no sinks or sources. Then
$H(V)=\set{(n,v)\in\N^k\times\Lambda}$ is the smallest hereditary
set containing $V$ and hence $\sat{V}=\Z^k\times\Lambda^0$ because
every vertex has a finite exhaustive set.

Conversely suppose $\KGalgCZLV$ is full and let $H(V)$ be the
hereditary set of $V$. Hence $H(V)\subset\N^k\times\Lambda^0$.
Then for every $(x,w)\in\SPgZL^0\setminus H(V)$ there must be a
path from $(x,w)$ to $H(V)$ and so there must exist $\lambda\in
w\Lambda$ such that $d(\lambda)+x\in\N^k$ for all $x\in\Z^k$.
Hence $w\Lambda^n$ is nonempty for all $w\in\Lambda^0$ and
$n\in\N^k$ and in particular $\Lambda^nv\ne\emptyset$ for any
$v\in H(V)$ and $n\in\N^k$. Therefore the set
$\set{v\in\Lambda^0:(x,v)\in H(V)}$ are non stranded vertices in
$\Lambda$ and every vertex in $\Lambda^0$ can trace a path to this
set. Hence $\Lambda$ is essentially saturated. Further for each
$(x,v)\in\SPgZL$ there exists a finite exhaustive subset of
$(x,v)\left(\SPgZL\right)H(V)$ for all $x\in\Z^k$ and hence there
is also exists a finite exhaustive set for each $v\in\Lambda^0$.
Thus $\Lambda$ is finitely exhaustive.
\end{proof}

\begin{cor} Given a row finite and essentially
saturated $k$-graph $\Lambda$ then the AF core $\KGalg^\gamma$ is
Morita equivalent to $\KGalg[\SPgZL]$.
\end{cor}

\begin{proof}
Follows from \propref{prop:AFcore} and \propref{prop:MEcore}.
\end{proof}

\newpage

\bibliographystyle{amsalpha}
\bibliography{Bib}

\providecommand{\bysame}{\leavevmode\hbox to3em{\hrulefill}\thinspace}
\providecommand{\MR}{\relax\ifhmode\unskip\space\fi MR }
\providecommand{\MRhref}[2]{%
  \href{http://www.ams.org/mathscinet-getitem?mr=#1}{#2}
}
\providecommand{\href}[2]{#2}
\begin{thebibliography}{KPRR}

\bibitem[APS]{APS}
S.~Allen, D.~Pask, and A.~Sims, \emph{A dual graph construction for higher rank
  graphs, and ${K}$-theory for finite 2-graphs}, Proc. Amer. Math. Soc. (2005),
  to appear [arXiv:math.OA/0402126].

\bibitem[Bat]{Ba}
T.~Bates, \emph{Applications of the gauge-invariant uniqueness theorem for
  graph algebras}, Bull. Austral. Math. Soc. \textbf{65} (2002), no.~1, 57--67.

\bibitem[BGR]{BGR}
L.~Brown, P.~Green, and M.~Rieffel, \emph{Stable isomorphism and strong
  {M}orita equivalence of ${C}^*$-algebras}, Pacific J. Math. \textbf{71}
  (1977), no.~2, 349--363.

\bibitem[BHRS]{BHRS}
T.~Bates, J.~H. Hong, I.~Raeburn, and W.~Szyma\'nski, \emph{The ideal structure
  of the ${C}^*$–algebras of infinite graphs}, Illinois J. Math. \textbf{46}
  (2002), 1159--1176.

\bibitem[BP]{BP}
T.~Bates and D.~Pask, \emph{Flow equivalence of graph algebras}, Ergodic Theory
  \& Dynam. Systems \textbf{24} (2004), no.~2, 367--382.

\bibitem[BPRS]{BPRS}
T.~Bates, D.~Pask, I.~Raeburn, and W.~Szyma\'nski, \emph{The ${C}^*$-algebras
  of row-finite graphs}, New York J. Math. \textbf{6} (2000), 307--324.

\bibitem[Bra]{Br}
O.~Bratteli, \emph{Inductive limits of finite dimensional ${C}^*$-algebras},
  Trans. Amer. Math. Soc. \textbf{17} (1972), 195--234.

\bibitem[Cri]{Cr}
T.~Crisp, \emph{Corners of graph algebras}, preprint, 2005
  [arXiv:math.OA/0503626].

\bibitem[DS]{DS}
D.~Drinen and N.~Sieben, \emph{${C}^*$-equivalences of graphs}, J. Operator
  Theory \textbf{45} (2001), no.~1, 209--229.

\bibitem[DT]{DT2}
D.~Drinen and M.~Tomforde, \emph{The ${C}^*$-algebras of arbitrary graphs},
  Rocky Mountain J. Math. \textbf{35} (2005), no.~1, 105--135.

\bibitem[KP1]{KP1}
A.~Kumjian and D.~Pask, \emph{${C}^*$-algebras of directed graphs and group
  actions}, Ergodic Theory \& Dynam. Systems \textbf{19} (1999), 1503--1519.

\bibitem[KP2]{KP2}
\bysame, \emph{Higher rank graph ${C}^*$-algebras}, New York J. Math.
  \textbf{6} (2000), 1--20.

\bibitem[KP3]{KP3}
\bysame, \emph{Actions of ${\Z}^k$ associated to higher rank graphs}, Ergod.
  Theory \& Dynam. Systems \textbf{23} (2003), no.~4, 1153--1172.

\bibitem[KPRR]{KPRR}
A.~Kumjian, D.~Pask, I.~Raeburn, and J.~Renault, \emph{Graphs, groupoids, and
  {C}untz-{K}rieger algebras}, J. Funct. Anal. \textbf{144} (1997), no.~2,
  505--541.

\bibitem[LM]{LM}
D.~Lind and B.~Marcus, \emph{An introduction to symbolic dynamics and coding},
  Cambridge University Press, 1995.

\bibitem[PQR]{PQR2}
D.~Pask, J.~Quigg, and I.~Raeburn, \emph{Coverings of $k$-graphs}, J. Algebra
  \textbf{289} (2005), no.~1, 161--191.

\bibitem[PR]{PR}
D.~Pask and I.~Raeburn, \emph{On the {K}-theory of {C}untz-{K}rieger algebras},
  Publ. RIMS Kyoto Univ. \textbf{32} (1996), 415--443.

\bibitem[PRY]{PRY}
D.~Pask, I.~Raeburn, and T.~Yeend, \emph{Actions of semigroups on directed
  graphs and their ${C}^*$-algebras}, J. Pure Appl. Algebra \textbf{159}
  (2001), no.~2-3, 297--313.

\bibitem[RS]{RS3}
G.~Robertson and T.~Steger, \emph{Asymptotic {K}-theory for groups acting on
  $\tilde{A}_2$ buildings}, Canadian Math. J. \textbf{8} (2001), 111--131.

\bibitem[RSY1]{RSY1}
I.~Raeburn, A.~Sims, and T.~Yeend, \emph{Higher-rank graphs and their
  ${C}^*$-algebras}, Proc. Edinburgh Math. Soc. \textbf{46} (2003), 99--115.

\bibitem[RSY2]{RSY2}
\bysame, \emph{The ${C}^*$-algebras of finitely aligned higher-rank graphs}, J.
  Func. Anal. \textbf{213} (2004), 206--240.

\bibitem[RW]{RW}
I.~Raeburn and D.~Williams, \emph{Morita equivalence and continuous-trace
  ${C}^*$-algebras}, Mathematical Surveys and Monographs, vol.~60, Amer. Math.
  Soc., 1998.

\bibitem[Sim]{Si2}
A.~Sims, \emph{Gauge-invariant ideals in the ${C}^*$-algebras of finitely
  aligned higher-rank graphs}, preprint, 2005[arXiv:math.OA/0406592].

\bibitem[Szy1]{Szy01}
W.~Szyma\'nski, \emph{Simplicity of {C}untz-{K}rieger algebras of infinite
  matrices}, Pacific J. Math. \textbf{199} (2001), no.~1, 249--256.

\bibitem[Szy2]{Szy02}
\bysame, \emph{The range of {$K$}-invariants for {$C^*$}-agebras of infinite
  graphs}, Ind. Univ. Math. J. \textbf{51} (2002), no.~1, 239--249.

\bibitem[Tom]{Tom1}
M.~Tomforde, \emph{Ext classes and embeddings for ${C}^*$-algebras of graphs
  with sinks}, New York J. Math. \textbf{7} (2001), 233--256.

\bibitem[Tyl]{Ty}
J.~Tyler, \emph{Every {AF}-algebra is {M}orita equivalent to a graph algebra},
  Bull. Austral. Math. Soc. \textbf{69} (2004), no.~2, 237--240.

\end{thebibliography}
\end{document}